\documentclass[twocolumn,showpacs,preprintnumbers,amsmath,amssymb]{revtex4}
\usepackage{graphicx}
\usepackage{dcolumn}
\usepackage{bm}
\usepackage{color}

\newcommand{\NaF}{\mathrm{Na}}
\newcommand{\KDR}{\mathrm{K}}
\newcommand{\KM}{\mathrm{M}}

\newcommand{\CaH}{\mathrm{Ca}}
\newcommand{\bd}{bifurcation diagram}
\newcommand{\Js}{J_{\star}}
\newcommand{\lc}{limit cycle}
\newcommand{\tc}{torus}

\newcommand{\JTB}{J_{\mathrm{TB}}}
\newcommand{\rd}{}

\begin{document}

\preprint{}

\title{New dynamics in cerebellar Purkinje cells:  torus canards}

\author{Mark A. Kramer$^{\dagger}$}
\author{Roger D. Traub}
\altaffiliation{Department of Physiology \& Pharmacology, SUNY Downstate Medical Center, Brooklyn, NY, 11203, USA}
\author{Nancy J. Kopell}
\altaffiliation{Department of Mathematics and Statistics, Boston University, Boston, MA, 02215, USA}


\begin{abstract}
We describe a transition from bursting to rapid spiking in a reduced mathematical model of a cerebellar Purkinje cell.  We perform a slow-fast analysis of the system and find that --- after a saddle node bifurcation of limit cycles --- the full model dynamics follow temporarily a repelling branch of limit cycles.  We propose that the system exhibits a dynamical phenomenon new to realistic, biophysical applications:  torus canards.
\end{abstract}


\maketitle
Bursting --- a repeated pattern of alternating quiescence and rapid spiking --- occurs in many neural systems, perhaps with functional implications \cite{Lisman:1997p5235, Izhikevich:2003p5531}.  Mathematical models developed to characterize bursting neural activity typically share common dynamical traits (e.g., excitability, slow-fast dynamics) and bifurcations \cite{Rinzel:1987fk, Izhikevich:2000p4238}.  In these models, the mechanisms that produce both periodic spiking and bursting activity are well understood \cite{HODGKIN:1952p5951, Izhikevich:2007fk}.  Yet the transition between these states often produces complicated dynamics (e.g., chaos, homoclinic bifurcations, blue sky catastrophes, period doubling cascades) more difficult to characterize \cite{TERMAN:1992p5295, Wang:1993p5394, Belykh:2000p5698, Cymbalyuk:2005, Shilnikov:2005p5302}.

In this letter, we describe a novel mechanism for the transition from bursting to spiking activity observed in a realistic, biophysical model of a cerebellar Purkinje cell.  We propose a reduction of this detailed model to study the transition and describe an intermediate state during which fast spiking activity is amplitude modulated by a slower rhythm.  In this intermediate state we observe a new type of dynamics in a continuous system that follows the attracting and repelling branches of {\it limit cycles} in the fast subsystem.  We compare these dynamics to traditional canard phenomena and propose that a new type of canard --- a torus canard --- occurs in the reduced model and provides a potential explanation of the detailed model activity.

\begin{figure}
\includegraphics[height=1.0in, width=3.4in]{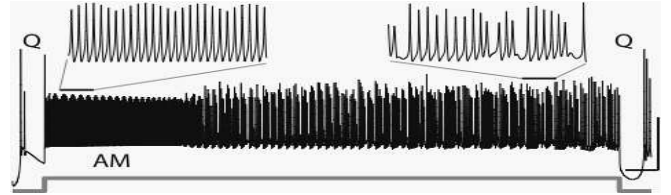}
\caption{\label{fig:traub} Simulation results of the soma compartment voltage (black)
from a detailed Purkinje cell model \cite{Middleton:2008fk}.  An injected current (gray) depolarizes the cell and produces rapid spiking modulated in amplitude (AM, upper left expanded trace) followed by more complicated activity (upper right expanded trace) between the quiescent intervals (Q) of bursts.  The horizontal and vertical scales (right) denote 50 ms and 50 mV, respectively.}
\end{figure}

The modeling and analysis were motivated by results observed in a detailed computational model of a cerebellar Purkinje cell \cite{Miyasho:2001p2686, Middleton:2008fk}.  The detailed model consists of $559$ compartments, each with $12$ types of ionic currents, resulting in over $6000$ dynamical variables.  We illustrate the results of a typical numerical simulation of this model in Figure \ref{fig:traub}.  Between the quiescent intervals of bursts (Q), we observe rapid spiking activity modulated in amplitude.  This modulation becomes more complicated as time progresses until the activity reenters the quiescent burst phase.  What dynamical and biophysical mechanisms produce this bursting activity interspersed with amplitude modulated (AM) spiking?

To answer this, we propose a reduced model of the detailed cerebellar Purkinje cell consisting of a single compartment with four ionic currents.  The voltage and current dynamics follow:
\begin{subequations}
\label{eq:whole}
\begin{eqnarray}
\!\!\!\!\!\!\dot{V} \!\! &=& \!\! -J  -g_{\KDR} n^4 (V\!+\!95)  -g_{\NaF}m_{0}[V]^3 h (V\!-\!50)   \label{eq:V}\\
\!\!\!\!\!\!&&\!\!\!\! -g_{\mathrm{L}}(V\!+\!70)\!-\!g_{\CaH} c^2 (V\!-\!125)\!-\!g_{\KM} M (V\!+\!95)  \nonumber\\
\!\!\!\!\!\!\dot{x} &=& ({x_{0}[V]-x})/{\tau_x[V]} \label{eq:X} .
\end{eqnarray}
\end{subequations}
Five currents affect the voltage in (\ref{eq:V}): a delayed rectifier potassium current ($g_{\KDR}$=10.0),
a transient inactivating sodium current ($g_{\NaF}$=125.0),
a leak current ($g_{\mathrm{L}}$=2.0),
a high-threshold noninactivating calcium current ($g_{\CaH}$=1.0),
and a muscarinic receptor suppressed potassium current (or M-current, $g_{\KM}$=0.75).  The dynamics of each gating variable follow (\ref{eq:X}) with $x$ replaced by $n, h, c,$ or $M$.  We implement the equilibrium function ($x_{0}[V]$) and time constant ($\tau_x[V]$) for each current from \cite{Traub:2003fk} and make the standard approximation of replacing the sodium activation variable with its equilibrium function ($m_{0}[V]$).  Of the five variables, the M-current evolves on a much slower time scale (at least ten times slower) and acts as the slow variable in this slow-fast system.  In what follows, we increase the excitation of the reduced model neuron by increasing the magnitude of parameter $J$ and compute numerical solutions and bifurcation diagrams for the system with XPPAUT and AUTO \cite{Ermentrout:2002p5541}.

\begin{figure}
\includegraphics[height=1.95in, width=3.4in]{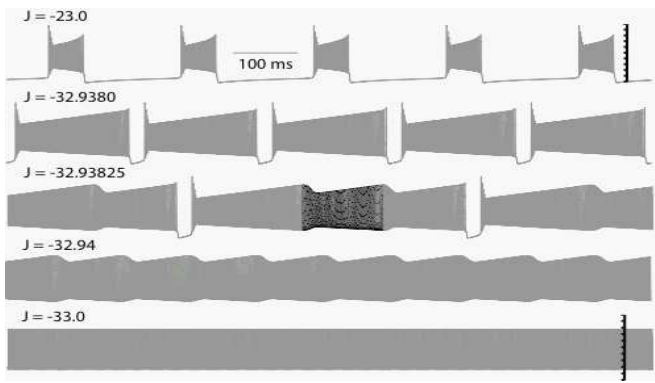}  
\caption{\label{fig:traces} The transition from bursting to AM spiking to rapid spiking in the reduced model as $J$ decreases.  We plot in gray the voltage activity of the neuron for five different values of $J$.  In the middle trace ($J=\Js$) we observe both bursting and AM spiking activity. We shade one complete AM cycle in black.  Further reductions in $J$ result in AM spiking and (unmodulated) spiking.  The vertical bars in the top and bottom traces denote a 100 mV range, extending from -60 mV to 40 mV, with 10 mV between tick marks.}
\end{figure}

We begin with a description of the voltage activity computed for decreasing values of the parameter $J$ to illustrate the transition from bursting to rapid spiking.  For $J > -22.5$, the dynamics approach a stable fixed point ($V \approx -54$ mV for $J=-22.5$, not shown.)  As we decrease $J$ through $-22.5$, bursts of activity emerge (Figure \ref{fig:traces}, top).  Within a burst the interval between, and amplitude of, the rapid spiking increase (from 1.6 ms to 2.0 ms, and 40 mV to 65 mV, respectively) after an initial transient.  The interburst intervals (lasting approximately 200 ms) are much longer than the intervals between spikes.  Decreasing $J$ further we find that the burst duty cycle increases, but that the interval between burst onsets remains approximately constant.  As we depolarize the neuron, more M-current must slowly accumulate to stop the bursting;  thus, the duration of spiking, compared to quiescence, increases.  Near $J=\Js=-32.93825$ the transition from bursting to rapid spiking begins and a new type of activity appears:  bursts interspersed with amplitude modulated (AM) fast spiking activity (Figure \ref{fig:traces}, middle).  The new activity increases the interval between bursts by integer multiples of 120 ms, the period of one complete AM cycle.  In Figure \ref{fig:traces}, one AM cycle separates the quiescent burst phases.  We find (but do not show) that the number of AM cycles between bursts appears unpredictable;  in simulations, we have observed between zero and 25 AM cycles between bursts.  Reducing the parameter further to $J=-32.94$ we find only AM spiking (and no bursting) activity.  For $J < -32.96$ only rapid spiking without amplitude modulation occurs and the transition from bursting to rapid spiking is complete.

What dynamical mechanisms govern the intermediate state between bursting and rapid spiking in the model?  To address this, we isolate the fast subsystem and examine its set of equilibria (i.e., the critical manifold) and periodic orbits.  In this slow-fast decomposition of (1), we fix $J$, treat the slow variable $M$ in (\ref{eq:whole}) as a parameter, and compute \bd s numerically in AUTO as implemented in XPPAUT \cite{Ermentrout:2002p5541}.  We then study the global dynamics of the full model on the \bd\ of the fast subsystem.  For $J= -23.0$ the dynamics exhibit well known behavior:  rapid spiking begins at a fold of fixed points and ends at a fold of limit cycles (i.e., fold / fold cycle bursting \cite{Izhikevich:2000p4238}).  But, as we decrease $J$ towards $\Js$, we find novel activity develops as we now describe.

\begin{figure}[t]
\includegraphics[height=2.0in, width=3.4in]{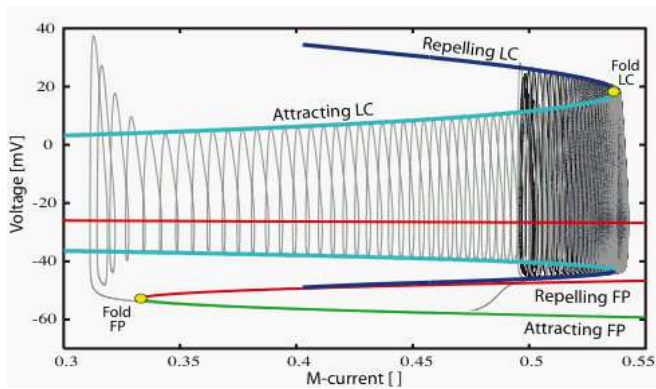}  
\caption{\label{fig:bif} Bifurcation diagram of the fast subsystem (thick and color) and dynamics of the full system (thin and gray/black) for $J=\Js$.  In the fast subsystem, the attracting and repelling fixed points and \lc s meet at folds (yellow circles) labeled Fold FP and Fold LC, respectively.  The attracting \lc s appear in a supercritical Hopf bifurcation (not shown).  At the fold of limit cycles the (slow) M-current changes direction and the complete dynamics follow the branch of repelling \lc s temporarily until returning to the curve of attracting fixed points (gray) or attracting \lc s (black).}
\end{figure}

In Figure \ref{fig:bif} we plot a bifurcation diagram (thick and color) for the fast subsystem and simulation results for the full system (thin and grayscale) with $J=\Js$.  In the full system, rapid spiking begins when the M-current decreases past the fold of fixed points in the fast subsystem;  at this fold or {\it knee}, the attracting and repelling fixed points merge and annihilate.  The voltage then increases rapidly, and the full dynamics approach the attracting curve of limit cycles in the fast subsystem.  With each spike, the slow M-current in the full system increases until the dynamics reach a fold of limit cycles in the fast subsystem.  At this fold, we expect spiking to cease consistent with the fold / fold cycle bursting observed for $J=-23.0$.  Instead we find that spiking continues as the dynamics of the full system follow the curve of limit cycles through the fold to the branch of {\it repelling} limit cycles.  The M-current decreases, and the full dynamics follow temporarily the repelling branch of limit cycles until returning to the branch of attracting fixed points (light gray) or limit cycles (black).  If the former, then the dynamics enter the quiescent phase of bursting and the M-current decreases.  If the latter, then an AM cycle occurs;  rapid spiking continues and the dynamics again approach the fold of limit cycles as the M-current increases.

\begin{figure}[t]
\includegraphics[height=3.36in, width=3.4in]{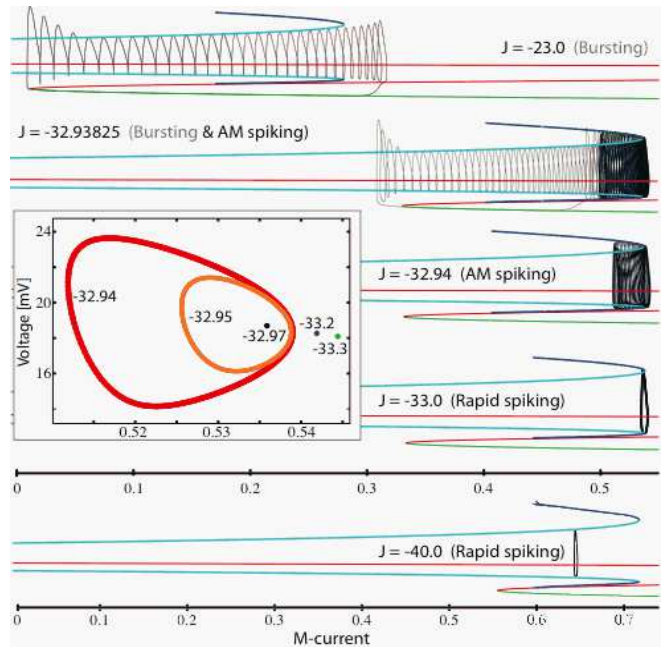}  
\caption{\label{fig:bif2} Bifurcation diagrams of the fast subsystem (thick and color) and dynamics of the full system (thin and gray/black) for five different values of $J$.  In the main figures, the vertical axes are identical;  the horizontal axes (identical for top four figures) indicate the value of $M$ and the axis is expanded for the bottom figure.  For $J=-23$, only fold / fold cycle bursting occurs (gray).  At $J=\Js$, the dynamics follow temporarily the branch of repelling \lc s, eventually returning to the branch of attracting fixed points (burst in gray) or the branch of attracting \lc s (AM cycle in black).  Reducing $J$ further results in AM spiking alone and rapid, unmodulated spiking.  {\rd Inset:  Poincare map (voltage versus M-current) of the full system for five different values of $J$.  Stable fixed points, corresponding to unmodulated rapid spiking, occur for three values of $J$.  For $J = \{-32.95,  -32.94\}$ stable invariant closed curves indicate amplitude modulation in the complete system.}}
\end{figure}

Decreasing $J$ past $\Js$ to $-32.94$ eliminates bursting in the full model dynamics and results in AM spiking alone.  During one complete AM cycle, the slow M-current in the full dynamics increases along the branch of attracting \lc s, passes through the fold of limit cycles, and decreases along the branch of repelling \lc s before returning to the branch of attracting \lc s (Figure \ref{fig:bif2}).  As $J$ decreases the extent of this slow modulation also decreases both in period (from approximately 0.11 s to 0.08 s) and in magnitude (the AM of the rapid spiking decreases from approximately 16 mV to 0.5 mV).  These reductions are suggested in Figure \ref{fig:bif2} and coincide with smaller excursions of the full dynamics from the fold point as $J$ decreases.  The transition from AM spiking to unmodulated spiking occurs at a supercritical torus bifurcation (negative first Lyapunov coefficient) near $J=\JTB=-32.96$ in the full system.  At this bifurcation, a stable torus and unstable limit cycle meet, and a stable limit cycle emerges.
{\rd The multipliers of the bifurcating \lc\ possess a complex conjugate pair whose moduli decrease through one as $J$ decreases through $\JTB$.  The resulting stable \lc s possess four multipliers of moduli less than one, and one multiplier of unit modulus corresponding to the fixed radius of the orbit.  We illustrate this transition in a Poincare map sampled at the apex of each spike in $V$ (Figure \ref{fig:bif2}).}  Decreasing $J$ further produces lower amplitude, faster oscillations that cease when $J < -150$ and the supercritical Hopf bifurcation occurs at $M > 1$, outside the physiological range.

We propose that the dynamics described above extend in new directions the classical canard phenomena observed in lower dimensional slow-fast systems.  The prototypical canard example consists of two variables evolving on different time scales (e.g., the van der Pol equation.)  For this two dimensional system, the critical manifold contains curves of attracting and repelling fixed points. In the full (2-D) system, the canards initiate at a Hopf bifurcation.  The subsequent periodic dynamics follow the curve of attracting fixed points until reaching a fold of fixed points in the fast subsystem.  Here, the slow variable reverses direction and the full dynamics follow temporarily the curve of repelling fixed points of the fast subsystem, eventually returning to a stable branch of fixed points \cite{DIENER:1984p5764}.

A related, but novel, phenomena appears to occur in the reduced, 5-dimensional model of the cerebellar Purkinje cell.  The bifurcation diagram in Figure \ref{fig:bif} illustrates the union of fixed points and limit cycles of the fast subsystem as $M$ is changed. During the active phase of the burst, the full model dynamics follow the curve of attracting limit cycles.  Each (fast) cycle within the burst increases the (slow) M-current until the global dynamics reach a fold in the \bd.  At this fold of {\it limit cycles} the average dynamics of the slow variable reverse direction (i.e., the dynamics of the M-current averaged over individual spikes in the fast subsystem change sign from positive to negative.)  The full dynamics then follow temporarily the curve of repelling {\it limit cycles}.  Consistent with classical canard phenomena, the parameter value ($J$) determines the length of time spent near the repelling branch and to which stable branch the dynamics return.  What differs here is that the canard initiates not after a Hopf bifurcation at a fold of fixed points, but instead after a torus bifurcation at a fold of \lc s.  We therefore label this phenomena a \tc\ canard.

The torus canard serves as an intermediary between the bursting and rapid spiking states.  This is often the case for canards associated with dramatic changes in dynamics resulting from small changes in a control parameter (e.g., a canard explosion).  In addition, complex (or chaotic) behavior often appears near the transition between different types of activity and may occur here (Figure 2, middle).  The sequence of bifurcations in this transitional region may be quite complicated and warrants further study.

The (slow) M-current and (fast) calcium current play complementary biophysical roles in the reduced model.  During the active phase of a burst, the M-current acts to hyperpolarize the cell and discourage spiking, while the calcium current acts to depolarize the cell and promote spiking.  When $J > \Js$, the hyperpolarization eventually wins, spiking stops, and the cell enters the quiescent phase of the burst.  The M-current is essential to this bursting activity \cite{Roopun:2006qy}.  For $J < \Js$, the depolarizing effect of the calcium current prevents the runaway hyperpolarization due to the M-current.  The quiescent phase of the bursts no longer occurs and we find instead only slow modulation of the fast spiking activity.  Decreasing $J$ further reduces the M-current dynamics and produces unmodulated fast spiking activity.  We note that, in the reduced model, blocking either the M or calcium current during AM spiking produces unmodulated rapid spiking.  These predictions were confirmed in the detailed model.

In this letter we described a novel mechanism that occurred during the transition from bursting to rapid spiking activity:  torus canards.  We identified this activity in a physiologically realistic computational model of a cerebellar Purkinje cell that motivated the simplified mathematical model.  By studying the slow-fast dynamics of the reduced system, we developed a better understanding of the physiological mechanisms (the M-current and calcium current) and dynamical mechanisms (the torus canard) that could produce the activity observed in the detailed model.

\vspace{-0.5mm}
Canards have been observed in other mathematical models of neural systems \cite{Guckenheimer:2000p5945,Rotstein:2003p3900,Moehlis:2006p4045,Rubin:2007p4333}.  However, in all of these systems, the canards occur along branches of attracting and repelling fixed points.  The dynamics presented here are unique in that the canards occur along branches of attracting and repelling limit cycles in a realistic, biophysical model.  In addition, these results were not limited to the reduced model;  similar dynamics were also observed in a detailed biophysical model.  Moreover, the slow modulation of the fast spiking activity appears to occur {\it in vitro} (see Figures 7A and 8A of \cite{Llinas:1980p5960}).  Additional recordings from cerebellar Purkinje cells could test experimentally the existence of torus canards and the role of the M-current in these dynamics.

\vspace{-0.5mm}
Our analysis focused on a computational, slow-fast decomposition of (1).  Although useful, this decomposition appears inadequate;  we note in particular that small changes in the middle three bifurcation diagrams of Figure \ref{fig:bif2} coincide with large changes in the full dynamics (namely, the transition from bursting to AM spiking to rapid spiking).  A better understanding of these dynamics will require a more sophisticated treatment \cite{Szmolyan:2001p6004}
involving, perhaps, dimensional reduction \cite{Medvedev:2006p5252} or an analysis of the canards in associated Poincare maps \cite{SHILNIKOV:2003p7899}.  The transition from bursting to spiking exhibits a mixed mode oscillation pattern (with bursts acting as the large amplitude oscillations and AM cycles as the small amplitude oscillations).  The system may therefore benefit from this type of analysis as well \cite{Brons:2006kx}.

\vspace{-0.5mm}
The authors thank Tasso J.\ Kaper and Horacio G.\ Rotstein for useful suggestions and discussions.  MAK is supported by NSF DMS-0602204.  RDT is supported by NIH R01NS04413.  NJK is supported by NSF DMS-0717670.

\end{document}